\documentclass[a4paper,twoside,english]{amsart}
\usepackage[active]{srcltx}
\usepackage{amsthm}
\usepackage{amssymb}
\usepackage{esint}
\usepackage{mathrsfs}

\makeatletter

\newcommand{\noun}[1]{\textsc{#1}}

\theoremstyle{plain}
\newtheorem{theorem}{\bf Theorem}[section]

\makeatother

\usepackage{babel}

\begin{document}

\title[A note on Boole polynomials with $q$-parameter]{A note on Boole polynomials with $q$-parameter}

\author[D. S. Kim,  Y. S. Jang,  T. Kim, and S.-H. Rim,]{ Dae San Kim, Yu Seon Jang, Taekyun Kim, and Seog-Hoon Rim,}
\begin{abstract}
Recently, Boole polynomials have been studied by Kim and Kim over the $p$-adic number field. In this paper, we consider a $q$-extension of Boole polynomials by using the fermionic $p$-adic integrals on $\mathbb{Z}_p$ and give some new identities related to those polynomials.
\end{abstract}

\keywords{}
\maketitle

\section{Introduction}
Let $p$ be a fixed odd prime number. Throughout this paper, $\mathbb{Z}_p$, $\mathbb{Q}_p$ and $\mathbb{C}_p$ will denote the ring of $p$-adic integers, the field of $p$-adic numbers and the completion of algebraic closure of $\mathbb{Q}_p$. The $p$-adic norm $|\cdot|_p$ is normalized as $|p|_p=1/p$. Let $C(\mathbb{Z}_p)$ be the space of continuous functions on $\mathbb{Z}_p$. For $f\in C(\mathbb{Z}_p)$, the fermionic $p$-adic integral on $\mathbb{Z}_p$ is defined by Kim to be

\begin{equation}
I(f)=\int_{\mathbb{Z}_p} f(x)d\mu_{-1}(x)= \lim_{N\rightarrow \infty}\sum_{x=0}^{p^N-1} f(x)(-1)^{x},\ (\textnormal{see}\ [11]).
\end{equation}
From (1), we have

\begin{equation}
I_{-1}(f_n)=2\sum_{x=0}^{n-1} (-1)^{n-1-a}f(a)+(-1)^nI_{-1}(f).
\end{equation}
Let us assume that $q$ is an indeterminate in $\mathbb{C}_p$ with $|1-q|_p<p^{-1/(p-1)}.$ The Stirling number of the first kind is defined by

\begin{equation}
(x)_n=\sum_{\ell=0}^n S_1(n,\ell)x^\ell, \ (n\geq0),
\end{equation}
and the Stirling number of the second kind is given by

\begin{equation}
x^n=\sum_{\ell=0}^{n}S_2(n, \ell)(x)_{\ell}, \ (n\geq0), \ (\textnormal{see}\ [9, 16]).
\end{equation}
The Boole polynomials are defined by the generating function to be

\begin{equation}
\sum_{n=0}^\infty Bl_n(x | \lambda)\frac{t^n}{n!}= \frac{1}{1+(1+t)^\lambda} (1+t)^x, \ (\textnormal{see}\ [9, 16]).
\end{equation}
In [9], Kim and Kim gave a Witt-type formula for $Bl_n(x | \lambda)$ over the $p$-adic number field as follows:

\begin{equation}
\int_{\mathbb{Z}_p} (x+\lambda y)_n d\mu_{-1}(y) = 2Bl_n(x | \lambda),\ (n\geq0),
\end{equation}
where $\lambda \in \mathbb{Z}_p$ and $(x)_n=x(x-1)\cdots(x-n+1)$.

Let us define the $q$-product of $x$ as follows:
\begin{equation}
(x)_{n,q}=x(x-q)(x-2q)\cdots(x-(n-1)q),\ (n\geq0).
\end{equation}
As is known, the Euler polynomials are defined by

\begin{equation}
\frac{2}{e^t+1} e^{xt}=\sum_{n=0}^\infty E_n (x) \frac{t^n}{n!}, \ (\textnormal{see}\ [1-20]).
\end{equation}
When $x=0, E_n=E_n(0)$ are called the Euler numbers. From (2), we have

\begin{equation}
\int_{\mathbb{Z}_p}e^{(x+y)t}d\mu_{-1}(y)=\frac{2}{e^t+1}e^{xt}=\sum_{n=0}^\infty E_n(x)\frac{t^n}{n!}.
\end{equation}

In this paper, we consider a $q$-extension of Boole polynomials by using the fermionic $p$-adic integrals on $\mathbb{Z}_p$ and give some new identities of those polynomials.

\section{Boole polynomials with $q$-parameter}

In this section, we assume that $t\in \mathbb{C}_p$ with $|t|_p<p^{-1/(p-1)}|q|_p$ and $\lambda \in \mathbb{Z}_p$. Now, we consider the Boole polynomials with $q$-parameter as follows:

\begin{equation}
Bl_{n,q}(x | \lambda)=\int_{\mathbb{Z}_p} \left(x+\lambda y \right)_{n,q} d\mu_{-1}(y),\ (n \geq 0).
\end{equation}
Thus, by (10), we get

\begin{equation}\begin{split}
Bl_{n, q}(x | \lambda)
=&\sum_{\ell=0}^n S_1(n, \ell) q^{n-\ell} \lambda^\ell \int_{\mathbb{Z}_p} \left(\frac{x}{\lambda}+y \right)^\ell d\mu_{-1}(y)\\
=&\sum_{\ell=0}^n S_1(n, \ell) q^{n-\ell} \lambda^\ell E_\ell \left( \frac{x}{\lambda}\right).
\end{split}\end{equation}
From (2) and (10), we can derive the generating function of $Bl_{n,q}(x | \lambda)$ as follows:

\begin{equation}\begin{split}
\sum_{n=0}^\infty Bl_{n,q}(x | \lambda)\frac{t^n}{n!}
=&\sum_{n=0}^\infty \int_{\mathbb{Z}_p}(x+\lambda y)_{n,q} d\mu_{-1}(y)\frac{t^n}{n!}\\
=&\sum_{n=0}^\infty q^n\int_{\mathbb{Z}_p} {\frac{x+\lambda y}{q} \choose n} d\mu_{-1}(y) t^n\\
=&\int_{\mathbb{Z}_p}(1+qt)^{\frac{x+\lambda y}{q}} d\mu_{-1}(y)\\
=&(1+qt)^{\frac{x}{q}}\left(\frac{2}{(1+qt)^{\lambda/q}+1}\right).
\end{split}\end{equation}
Therefore, by (12), we obtain the following theorem.

\begin{theorem} Let $F(t, x|\lambda)=\sum_{n=0}^\infty Bl_{n, q}(x | \lambda) \frac{t^n}{n!}.$ Then we have
\begin{equation*}
F(t, x | \lambda)=\left(\frac{2}{(1+qt)^{\lambda/q}+1} \right)(1+qt)^{x/q}.
\end{equation*}
\end{theorem}
By replacing $t$ by $(e^t-1)/q$ in (12), we get
\begin{equation}\begin{split}
\sum_{n=0}^\infty Bl_{n,q}(x | \lambda)q^{-n}\frac{(e^t-1)^n}{n!}
=&\frac{2}{e^{\lambda t/q }+1}e^{\frac{x}{q}t}\\
=&\sum_{n=0}^\infty E_n\left(\frac{x}{\lambda} \right) \left(\frac{\lambda}{q}\right)^n \frac{t^n}{n!}
\end{split}\end{equation}
and
\begin{equation}\begin{split}
\sum_{n=0}^\infty Bl_{n,q} (x|\lambda) q^{-n} \frac{(e^t-1)^n}{n!}
=&\sum_{n=0}^\infty Bl_{n,q} (x|\lambda) q^{-n} \sum_{m=n}^\infty S_2(m,n)\frac{t^m}{m!}\\
=&\sum_{m=0}^\infty \left\{ \sum_{n=0}^m Bl_{n,q}(x|\lambda) \frac{S_2(m,n)}{q^n}\right\} \frac{t^m}{m!}.
\end{split}\end{equation}
Therefore, by (13) and (14), we obtain the following theorem.

\begin{theorem} For $m\geq0$, we have
\begin{equation*}
\lambda^m E_m \left(\frac{x}{\lambda}\right) = \sum_{n=0}^m Bl_{n,q} (x|\lambda) q^{m-n} S_2(m,n),
\end{equation*}
and
\begin{equation*}
Bl_{m,q}(x|\lambda)=\sum_{\ell=0}^m S_1(m,\ell)q^{m-\ell}\lambda^\ell E_\ell \left(\frac{x}{\lambda} \right).
\end{equation*}
\end{theorem}
Note that $\lim_{q \rightarrow 1}Bl_{n,q}(x | \lambda)=2Bl_{n} (x | \lambda),\ (n \geq 0)$. When $x=0, Bl_{n,q} (\lambda)=Bl_{n,q}(0|\lambda)$ are called the $q$-Boole numbers.

Now, we consider the $q$-Boole polynomials of the second kind as follows:

\begin{equation}
\widehat{Bl}_{n,q}(x | \lambda)=\int_{\mathbb{Z}_p}(-\lambda y+x)_{n,q} d\mu_{-1} (y), \ (n \geq 0).
\end{equation}
Thus, by (15), we get

\begin{equation}\begin{split}
\widehat {Bl}_{n,q}(x | \lambda)
=& q^n \int_{\mathbb{Z}_p} \left( \frac{-\lambda y +x}{q}\right)_n d\mu_{-1}(y)\\
=& q^n \int_{\mathbb{Z}_p} \sum_{\ell=0}^n  \frac{\lambda^\ell S_1(n, \ell)}{q^\ell} (-1)^\ell \left(y-\frac{x}{\lambda} \right)^\ell d\mu_{-1}(y)\\
=& \sum_{\ell=0}^n S_1(n,\ell)q^{n-\ell} \lambda^\ell (-1)^\ell E_{\ell}\left( -\frac{x}{\lambda}\right).
\end{split}\end{equation}
From (8), we have

\begin{equation}\begin{split}
\sum_{n=0}^\infty E_n\left( -\frac{x}{\lambda}\right)\frac{t^n}{n!}
=&\frac{2}{e^t+1}e^{\left( -\frac{x}{\lambda}\right)t}\\
=&\frac{2}{1+e^{-t}}e^{-\left(1+ \frac{x}{\lambda}\right)t}\\
=&\sum_{n=0}^\infty (-1)^n E_n \left(1+\frac{x}{\lambda} \right)\frac{t^n}{n!}.
\end{split}\end{equation}
By (16) and (17), we get

\begin{equation}\begin{split}
\widehat{Bl}_{n,q} (x | \lambda)
=&\sum_{\ell=0}^n\lambda^\ell |S_1(n,\ell)|q^{n-\ell} E_{\ell}\left(1+\frac{x}{\lambda} \right).
\end{split}\end{equation}
From (15), we can derive the generating function of the Boole polynomials of the second kind as follows:
\begin{equation}\begin{split}
\sum_{n=0}^\infty \widehat{Bl}_{n,q} (x | \lambda)\frac{t^n}{n!}
=& \sum_{n=0}^\infty q^n \int_{\mathbb{Z}_p} {\frac{-\lambda y+x}{q} \choose n}  d\mu_{-1}(y) t^n\\
=& \int_{\mathbb{Z}_p} (1+qt)^{\frac{-\lambda y+x}{q}} d\mu_{-1}(y)\\
=& (1+qt)^{\frac{x}{q}}\int_{\mathbb{Z}_p} (1+qt)^{-\frac{\lambda y}{q}}d\mu_{-1}(y) \\
=& (1+qt)^{\frac{x+\lambda}{q}} \frac{2}{(1+qt)^{\lambda/q}+1}.
\end{split}\end{equation}
By replacing $t$ by $(e^t-1)/q$ in (19), we get
\begin{equation}\begin{split}
\sum_{n=0}^\infty \widehat{Bl}_{n,q} (x | \lambda) \frac{1}{q^n} \frac{(e^t-1)^n}{n!}
=& e^{\frac{1}{q}(x+\lambda)t} \frac{2}{e^{\frac{\lambda}{q}t}+1}\\
=& \sum_{n=0}^\infty E_n \left( \frac{x+\lambda}{\lambda}\right) \frac{\lambda^n}{q^n}\frac{t^n}{n!}
\end{split}\end{equation}
and
\begin{equation}\begin{split}
\sum_{n=0}^\infty \widehat{Bl}_{n,q} (x | \lambda) q^{-n} \frac{(e^t-1)^n}{n!}
=& \sum_{n=0}^\infty \widehat{Bl}_{n,q} (x | \lambda) q^{-n} \sum_{m=n}^\infty S_2 (m,n) \frac{t^m}{m!}\\
=& \sum_{m=0}^\infty \left(\sum_{n=0}^m \widehat{Bl}_{n,q} (x | \lambda) \frac{S_2(m,n)}{q^n} \right) \frac{t^m}{m!}.
\end{split}\end{equation}
Therefore, by (18), (20) and (21), we obtain the following theorem.
\begin{theorem} For $m \geq 0,$ we have
\begin{equation*}
\sum_{n=0}^m \widehat{Bl}_{n,q} (x | \lambda) q^{m-n} S_2 (m, n) = \lambda^m E_m \left(1+\frac{x}{\lambda}\right)
\end{equation*}
and
\begin{equation*}
\widehat{Bl}_{m,q} (x | \lambda) = \sum_{\ell=0}^m S_1(m, \ell) q^{m-\ell} \lambda^{\ell} E_{\ell} \left(1+\frac{x}{\lambda} \right).
\end{equation*}
\end{theorem}
For $\alpha\in \mathbb{N}$, let us consider $q$-Boole polynomials of the first kind with order $\alpha$ as follows:
\begin{equation}\begin{split}
Bl_{n,q}^{(\alpha)}(x | \lambda)
=& \int_{\mathbb{Z}_p}\cdots \int_{\mathbb{Z}_p} \left(\lambda x_1+\cdots+\lambda x_\alpha+x \right)_{n,q} d\mu_{-1}(x_1) \cdots d\mu_{-1}(x_\alpha)\\
=& q^n \int_{\mathbb{Z}_p}\cdots \int_{\mathbb{Z}_p} \left(\frac{\lambda x_1+\cdots+\lambda x_\alpha+x}{q} \right)_n d\mu_{-1}(x_1) \cdots d\mu_{-1}(x_\alpha)\\
=& q^n \sum_{\ell=0}^n S_1(n,\ell) \frac{1}{q^\ell}  \int_{\mathbb{Z}_p}\cdots \int_{\mathbb{Z}_p} \left(\lambda x_1+\cdots+\lambda x_\alpha+x \right)^\ell  d\mu_{-1}(x_1) \cdots d\mu_{-1}(x_\alpha)\\
=& \sum_{\ell=0}^n S_1 (n, \ell) q^{n-\ell} \lambda^\ell E_\ell ^{(\alpha)}\left( \frac{x}{\lambda}\right),
\end{split}\end{equation}
where $E_n^{(\alpha)}(x)$ are the Euler polynomials of order $\alpha$ which are defined by
\begin{equation*}
\left( \frac{2}{e^t+1}\right)^\alpha e^{xt} =\sum_{n=0}^\infty E_n ^{(\alpha)}(x)\frac{t^n}{n!}.
\end{equation*}
From (22), we note that the generating function of $Bl_{n,q}^{(\alpha)} (x | \lambda)$ are given by
\begin{equation}\begin{split}
\sum_{n=0}^\infty Bl_{n,q}^{(\alpha)} (x | \lambda) \frac{t^n}{n!}
=& \sum_{n=0}^\infty q^n \int_{\mathbb{Z}_p} \cdots \int_{\mathbb{Z}_p} {\frac{\lambda x_1+\cdots+\lambda x_\alpha+x}{q} \choose n } d\mu_{-1}(x_1) \cdots d\mu_{-1}(x_\alpha)t^{n}\\
=& \int_{\mathbb{Z}_p}\cdots \int_{\mathbb{Z}_p} (1+qt)^{\frac{\lambda x_1+ \cdots +\lambda x_\alpha +x}{q}} d\mu_{-1}(x_1) \cdots d\mu_{-1}(x_\alpha)\\
=&(1+qt)^{\frac{x}{q}} \left( \frac{2}{(1+qt)^{\frac{\lambda}{q}}+1}\right)^\alpha.
\end{split}\end{equation}
By replacing $t$ by $(e^t-1)/q$, we get
\begin{equation}\begin{split}
\sum_{n=0}^\infty Bl_{n,q}^{(\alpha)} (x | \lambda) \frac{1}{q^n} \frac{(e^t-1)^n}{n!}
=& e^{\frac{x}{q}t} \left(\frac{2}{e^{\frac{\lambda}{q}t}+1} \right)^\alpha\\
=& \sum_{n=0}^\infty E_n ^{(\alpha)} \left( \frac{x}{\lambda}\right) \frac{\lambda^n}{q^n} t^n
\end{split}\end{equation}
and
\begin{equation}\begin{split}
\sum_{n=0}^\infty Bl_{n,q} ^{(\alpha)} (x | \lambda) \frac{1}{q^n} \frac{(e^t-1)^n}{n!}
=& \sum_{n=0}^\infty Bl_{n,q}^{(\alpha)} (x | \lambda) \frac{1}{q^n} \sum_{m=n}^\infty S_2(m,n) \frac{t^m}{m!}\\
=& \sum_{m=0}^\infty \left(\sum_{n=0}^m Bl_{n,q}^{(\alpha)} (x | \lambda) \frac{S_2(m,n)}{q^n} \right)\frac{t^m}{m!}.
\end{split}\end{equation}
Therefore, by (24) and (25), we obtain the following theorem.
\begin{theorem} For $m\geq0$, we have
\begin{equation*}
\lambda^m E_m^{(\alpha)} \left( \frac{x}{\lambda}\right)=\sum_{n=0}^m q^{m-n} Bl_{n,q}^{(\alpha)} (x | \lambda) S_2(m,n)
\end{equation*}
and
\begin{equation*}
Bl_{m,q}^{(\alpha)} (x | \lambda)= \sum_{\ell=0}^m S_1 (m, \ell) q^{m-\ell} \lambda^\ell E_\ell ^{(\alpha)} \left(\frac{x}{\lambda} \right).
\end{equation*}
\end{theorem}
We consider the $q$-Boole polynomials of the second kind with order $\alpha$ as follows:
\begin{equation}
\widehat{Bl}_{n,q}^{(\alpha)} (x | \lambda)=\int_{\mathbb{Z}_p} \cdots \int_{\mathbb{Z}_p} \left( -\lambda x_1  - \cdots -\lambda x_\alpha +x\right)_{n,q}d \mu_{-1} (x_1) \cdots d\mu_{-1} (x_\alpha).
\end{equation}
Then, by (26), we get
\begin{equation}\begin{split}
\widehat{Bl}_{n,q} ^{(\alpha)} (x | \lambda)
=& q^n\sum_{\ell=0}^n S_1(n, \ell)(-1)^\ell \left( \frac{\lambda}{q} \right)^\ell \int_{\mathbb{Z}_p}\cdots \int_{\mathbb{Z}_p} \left(x_1 + \cdots +x_\alpha -\frac{x}{\lambda} \right)^\ell d\mu_{-1} (x_1) \cdots d\mu_{-1}(x_\alpha)\\
=& \sum_{\ell=0}^n S_1 (n,\ell)(-1)^\ell q^{n-\ell}\lambda^\ell E_{\ell} ^{(\alpha)} \left( -\frac{x}{\lambda}\right).
\end{split}\end{equation}
From the definition of the higher-order Euler polynomials, we note that
\begin{equation}\begin{split}
\sum_{n=0}^\infty E_n ^{(\alpha)} \left(-\frac{x}{\lambda} \right) \frac{t^n}{n!}
=& \left( \frac{2}{e^t+1}\right)^\alpha e^{- \frac{x}{\lambda}t}\\
=& \left(\frac{2}{1+e^{-t}} \right)^\alpha e^{- \left(\frac{x}{\lambda}+\alpha \right)t} \\
=& \sum_{n=0}^\infty (-1)^n E_n^{(\alpha)} \left(\frac{x}{\lambda}+\alpha \right) \frac{t^n}{n!}.
\end{split}\end{equation}
Thus, by (27) and (28), we get
\begin{equation}\begin{split}
\widehat{Bl}_{n,q}^{(\alpha)} (x | \lambda)
= &\sum_{\ell=0}^n S_1 (n, \ell) q^{n-\ell} \lambda^\ell E_\ell ^{(\alpha)} \left( \frac{x}{\lambda}+\alpha\right).
\end{split}\end{equation}
From (26), we have

\begin{equation}\begin{split}
\sum_{n=0}^\infty \widehat{Bl}_{n,q}^{(\alpha)} (x | \lambda) \frac{t^n}{n!}
=& \sum_{n=0}^{\infty}q^{n} \int_{\mathbb{Z}_p}\cdots \int_{\mathbb{Z}_p} {\frac{-\lambda x_1+\cdots-\lambda x_\alpha +x}{q} \choose n} d\mu_{-1}(x_1) \cdots d\mu_{-1}(x_\alpha)t^{n}\\
=& \int_{\mathbb{Z}_p}\cdots\int_{\mathbb{Z}_p} (1+qt)^{\frac{x-\lambda x_1 -\cdots- \lambda x_\alpha}{q}} d\mu_{-1}(x_1) \cdots d\mu_{-1}(x_\alpha)\\
= & (1+qt)^{\frac{x+\alpha}{q}} \left( \frac{2}{(1+qt)^{\frac{\lambda}{q}}+1}\right)^\alpha.
\end{split}\end{equation}
By replacing $t$ by $(e^t-1)/q$, we get
\begin{equation}\begin{split}
\sum_{n=0}^\infty \widehat{Bl}_{n,q}^{(\alpha)} (x | \lambda) \frac{(e^t-1)^n}{n!q^n}
=& e^{\frac{x+\alpha}{q}t} \left( \frac{2}{e^{\frac{\lambda}{q}t}+1}\right)^\alpha\\
=& \sum_{n=0}^\infty E_n ^{(\alpha)} \left( \frac{x+\alpha}{\lambda}\right) \left(\frac{\lambda}{q} \right)^n \frac{t^n}{n!}.
\end{split}\end{equation}
and
\begin{equation}\begin{split}
\sum_{n=0}^\infty \widehat{Bl}_{n,q}^{(\alpha)} (x | \lambda) \frac{1}{q^n} \frac{1}{n!} (e^t-1)^n
=& \sum_{n=0}^\infty \widehat{Bl}_{n,q}^{(\alpha)} (x | \lambda) \frac{1}{q^n} \sum_{m=n}^\infty S_2(m,n) \frac{t^m}{m!}\\
=& \sum_{m=0}^\infty \left( \sum_{n=0}^m \widehat{Bl}_{n,q}^{(\alpha)}(x | \lambda) \frac{S_2(m,n)}{q^n}\right) \frac{t^m}{m!}.
\end{split}\end{equation}
Therefore, by (31) and (32), we obtain the following theorem.

\begin{theorem} For $m\geq0$, we have
\begin{equation}
\lambda^m E_m^{(\alpha)} \left(\frac{x+\alpha}{\lambda} \right) = \sum_{n=0}^m q^{m-n} \widehat{Bl}_{n,q}^{(\alpha)} (x | \lambda) S_2 (m,n)
\end{equation}
and
\begin{equation}
\widehat{Bl}_{m,q} ^{(\alpha)} (x | \lambda) = \sum_{\ell=0}^m S_1(m, \ell) q^{m-\ell} \lambda^\ell E_\ell^{(\alpha)} \left(\frac{x+\alpha}{\lambda}\right).
\end{equation}
\end{theorem}

{\hskip -1pc \bf Remark.} When $x=0$, $\widehat{Bl}_{n,q}(\lambda)=\widehat{Bl}_{n,q}(0 | \lambda)$ are called the $q$-Boole numbers of the second kind.

Now, we observe that
\begin{equation}\begin{split}
\frac{\widehat{Bl}_{n,q}(\lambda)}{n!}
=& \frac{1}{n!} \int_{\mathbb{Z}_p} (-\lambda y)_{n,q} d\mu_{-1}(y)\\
=& q^n \int_{\mathbb{Z}_p} {\frac{-\lambda y}{q} \choose n} d\mu_{-1}(y)\\
=& (-1)^n q^n \int_{\mathbb{Z}_p} {\frac{\lambda y}{q}+n-1 \choose n} d\mu_{-1}(y)\\
=& (-1)^n q^n \sum_{\ell=0}^n {n-1 \choose \ell-1} \frac{1}{\ell!}\int_{\mathbb{Z}_p} \left(\frac{\lambda y}{q}\right)_\ell d\mu_{-1}(y)\\
=& (-1)^n q^n \sum_{\ell=0}^n {n-1 \choose \ell-1} \frac{Bl_{\ell, q}(\lambda)}{\ell!q^{\ell}}.
\end{split}\end{equation}
Therefore, by (35), we obtain the following theorem.

\begin{theorem} For $n\geq0$, we have
\begin{equation}
\frac{(-1)^n}{q^n} \frac{\widehat{Bl}_{n,q}(\lambda)}{n!} =\sum_{\ell=0}^n {n-1 \choose \ell-1} \frac{Bl_{\ell,q}(\lambda)}{\ell!q^{\ell}}.
\end{equation}
\end{theorem}

\bigskip
\medskip
\bigskip
\medskip

\noindent \noun{D. S. Kim\\
Department of Mathematics\\
Sogang Uiversity\\
Seoul, Republic of Korea}\\
\textnormal{e-mail: dskim@sogang.ac.kr}

\vskip 1pc

\noindent \noun{Y. S. Jang\\
Department of Applied Mathematics\\
Kangnam University\\
Yongin 446-702, Republic of Korea}\\
\textnormal{e-mail: ysjang@kangnam.ac.kr}

\vskip 1pc

\noindent \noun{T. Kim\\
 Department of Mathematics\\
Kwangwoon University\\
Seoul 139-701, Republic of Korea}\\
\textnormal{e-mail: tkkim@kw.ac.kr}

\vskip 1pc

\noindent \noun{S.-H. Rim\\
Department of Mathematics Education\\
Kyungpook National University\\
Taegu 702-701, Republic of Korea}\\
\textnormal{e-mail: shrim@knu.ac.kr}


\bigskip
\medskip

\begin{thebibliography}{99}
\bibitem{ref-1} S. Araci and M. Acikgoz, A note on the Frobenius-Euler numbers and polynomials associated with Bernstein polynomials, \emph{Adv. Stud. Contemp. Math.} {\bf 22} (2012), no. 3, 399--406.
\bibitem{ref-2} A. Bayad and J. Chikhi, Apostol-Euler polynomials and asymptotics for negative binomial reciprocals, \emph{Adv. Stud. Contemp. Math.} {\bf 24} (2014), no. 1, 33--37.
\bibitem{ref-3} M. Can, M. Cenkci, V. Kurt, and Y. Simsek, Twisted Dedekind type sums associated with Barnes' type multiple Frobenius-Euler $\ell$-functions, \emph{Adv. Stud. Contemp. Math.} {\bf 18} (2009), no. 2, 135--160.
\bibitem{ref-4} J. Choi, D. S. Kim, T. Kim, and Y. H. Kim, Some arithmetic identities on Bernoulli and Euler numbers arising from the $p$-adic integrals on $\mathbb{Z}_p$, \emph{Adv. Stud. Contemp. Math.} {\bf 22} (2012), no. 2, 239--247.
\bibitem{ref-5} D. Ding and J. Yang, Some identities related to the Apostol-Euler and Apostol-Bernoulli polynomials, \emph{Adv. Stud. Contemp. Math.} {\bf 20} (2010), no. 1, 7--21.
\bibitem{ref-6} D. V. Dolgy, T. Kim, B. Lee, and C. S. Ryoo, On the $q$-analogue of Euler measure with weight $\alpha$, \emph{Adv. Stud. Contemp. Math.} {\bf 21} (2011), no. 4, 429--435.
\bibitem{ref-7} S. Gaboury, R. Tremblay, and B.-J. Fug$\grave{e}$re, Some explicit foumulas for certain new classes of Bernoulli, Euler and Genocchi polynomials, \emph{Proc. Jangjeon Math. Soc.} {\bf 17} (2014), no. 1, 115--123.
\bibitem{ref-8} J. H. Jeong, J.-H. Jin, J.-W. Park, and S.-H. Rim, On the twisted weak $q$-Euler numbers and polynomials with weight $0$, \emph{Proc. Jangjeon Math. Soc.} {\bf 16} (2013), no. 2, 157--163.
\bibitem{ref-9} D. S. Kim and T. Kim, A note on Boole polynomials, \emph{Integral Transforms Spec. Funct.} {\bf 25} (2014), no. 8, 627--633.
\bibitem{ref-10} T. Kim, D. S. Kim, D. V. Dolgy, and S.-H. Rim, Some identities on the Euler numbers arising from Euler basis polynomials, \emph{Ars Combin.} {\bf 109} (2013), 433--446.
\bibitem{ref-11} T. Kim, Symmetry of power sum polynomials and multivariate fermionic $p$-adic invariant integral on $\mathbb{Z}_p$, \emph{Russ. J. Math. Phys.} {\bf 16} (2009), no. 1, 93--96.
\bibitem{ref-12} T. Kim, Some identities on the $q$-Euler polynomials of higher order and $q$-Stirling numbers by the fermionic $p$-adic integral on $\mathbb{Z}_p$, \emph{Russ. J. Math. Phys.} {\bf 16} (2009), no. 4, 484--491.
\bibitem{ref-13} T. Kim, A study on the $q$-Euler numbers and the fermionic $q$-integral of the product of several type $q$-Bernstein polynomials on $\mathbb{Z}_p$,  \emph{Adv. Stud. Contemp. Math.} {\bf 23} (2013), no. 1, 5--11.
\bibitem{ref-14} Q.-M. Luo, $q$-analogues of some results for the Apostol-Euler polynomials,  \emph{Adv. Stud. Contemp. Math.} {\bf 20} (2010), no. 1, 103--113.
\bibitem{ref-15} H. Ozden, I. N. Cangul, and Y. Simsek, Remarks on $q$-Bernoulli numbers associated with Daehee numbers, \emph{Adv. Stud. Contemp. Math.} {\bf 18} (2009), no. 1, 41--48.
\bibitem{ref-16} S. Roman, The umbral calculus, Academic Press, Inc., New York, 1984.
\bibitem{ref-17} C. S. Ryoo, H. Song and R. P. Agarwal, On the roots of the $q$-analogue of Euler-Barnes' polynomials,  \emph{Adv. Stud. Contemp. Math.} {\bf 9} (2004), no. 2, 153--163.
\bibitem{ref-18} Y. Simsek, O. Yurekli and V. Kurt, On interpolation functions of the twisted generalized Frobenius-Euler numbers,  \emph{Adv. Stud. Contemp. Math.} {\bf 15} (2007), no. 2, 187--194.
\bibitem{ref-19} E. Sen, Theorems on Apostol-Euler polynomials of higher order arising from Euler basis, \emph{Adv. Stud. Contemp. Math.} {\bf 23} (2013), no. 2, 337--345.
\bibitem{ref-20} Z. Zhang and H. Yang, Some closed formulas for generalized Bernoulli-Euler numbers and polynomials,  \emph{Proc. Jangjeon Math. Soc.} {\bf 11} (2008), no. 2, 191--198.

\end{thebibliography}
\end{document}